\documentclass{amsart}
\usepackage{amssymb}
\usepackage[all]{xy}
\newtheorem{theorem}{Theorem}[section]

\newtheorem{prop}[theorem]{Proposition}
\newtheorem{cor}[theorem]{Corollary}
\newtheorem{conjecture}[theorem]{Conjecture}
\theoremstyle{definition}
\newtheorem{definition}[theorem]{Definition}
\newtheorem{example}[theorem]{Example}

\newtheorem{question}[theorem]{Question}

\theoremstyle{remark}

\numberwithin{equation}{subsection}

\def\aa{\mathbb{A}}
\def\cc{\mathbb{C}}
\def\ff{\mathbb{F}}
\def\GG{\mathbb{G}}

\def\pp{\mathbb{P}}
\def\rr{\mathbb{R}}
\def\qq{\mathbb{Q}}
\def\zz{\mathbb{Z}}

\def\frF{\mathfrak{F}}

\def\deg{\mathrm{deg}}
\def\disc{\mathrm{disc}}
\def\dim{\mathrm{dim}}
\def\Frob{\mathrm{Frob}}
\def\Gal{\mathrm{Gal}}

\def\Id{\mathrm{Id}}
\def\inf{\mathrm{inf}}

\def\mod{\mathrm{mod}}
\def\ord{\mathrm{ord}}
\def\Sym{\mathrm{Sym}}

\begin{document}

\title{Lectures on Zeta Functions over Finite Fields}
\author{Daqing Wan}
\address{Department of Mathematics, University of California,
Irvine, Ca92697-3875}
\address{Email: dwan@math.uci.edu}
\date{}




\maketitle
\setcounter{tocdepth}{1}
\tableofcontents
\section{Introduction}
These are the notes from the summer school in G\"ottingen sponsored
by NATO Advanced Study Institute on Higher-Dimensional Geometry over
Finite Fields that took place in 2007. The aim was to give a short
introduction on zeta functions over finite fields, focusing on
moment zeta functions and zeta functions of affine toric
hypersurfaces. Along the way, both concrete examples and open
problems are presented to illustrate the general theory. For
simplicity, we have kept the original lecture style of the notes. It
is a pleasure to thank Phong Le for taking the notes and for his
help in typing up the notes.


\section{Zeta functions over Finite fields}
\subsection*{Definitions and Examples}
Let $p$ be a prime, $q=p^a$ and $\ff_q$ be the finite field of $q$ elements.  For the affine line $\aa^1$,
we have $\aa^1(\ff_q) = \ff_q$ and $\# \aa^1(\ff_q)=q$.

Fix an algebraic closure $\overline{\ff_q}$.  $\Frob_q:\overline{\ff_q} \mapsto \overline{\ff_q}$,
defined by $\Frob_q (x)= x^q$.  For $k \in \zz_{>0}$,
$$\ff_{q^k} = \mathrm{Fix}\left(\Frob_{q}^k |
\overline{\ff_q}\right), \ \ \aa^1(\overline{\ff_q}) =
\overline{\ff_q} = \bigcup_{k=1}^{\infty} \ff_{q^k}.$$

Given a geometric point $x \in \overline{\ff_q}$, the orbit $\{x ,
x^q,\ldots, x^{q^{\deg(x)-1}}\}$ of $x$ under $\Frob_q$ is called
the closed point of $\aa^1$ containing $x$.  The length of the orbit
is called the degree of the closed point. We may correspond this
uniquely to the monic irreducible polynomial $ (t-x)(t-x^q) \ldots
(t-x^{q^{\deg(x)-1}})$. Let $|\aa^1|$ denote the set of closed
points of $\aa^1$ over $\ff_q$. Similarly, let $|\aa^1|_k$ denote
the set of closed points of $\aa^1$ of degree $k$. Hence
$$|\aa^1| = \coprod_{k=1}^{\infty} |\aa^1|_k.$$

\begin{example}
The zeta function of $\aa^1$ over $\ff_q$ is
$$
\begin{array}{lll}
Z(\aa^1,T) & = & \exp\left( \sum_{k=1}^{\infty} \frac{T^k}{k} \#\aa^1(\ff_{q^k}) \right)\\
&= &\exp\left( \sum_{k=1}^{\infty} \frac{T^k}{k} q^k \right)\\
&= &\frac{1}{1-qT} \in \qq(T).
\end{array}
$$
The reciprocal pole is a Weil $q$-number.
There is also a product decomposition
$$Z(\aa^1,T) = \prod_{k=1}^{\infty} \frac{1}{(1-T^k)^{\#|\aa^1|_k}}.$$
\end{example}

More generally, let $X$ be quasi-projective over $\ff_q$, or a
scheme of finite type over $\ff_q$. By birational equivalence and
induction, one can often (but not always) assume that $X$ is a
hypersurface $\{f(x_1,\ldots,x_n)=0 | x_i \in \overline{\ff_q}\}$.
Consider the Frobenius action on $X(\overline{\ff_q})$.  Let $|X|$
be the set of all closed points of $X$ and $|X|_k$ be the set of
closed points on $X$ of degree $k$. As in the previous case, we
have
$$X(\overline{\ff_q}) = \bigcup_{k=1}^{\infty} X(\ff_{q^k}), \ \ |X| = \coprod_{k=1}^{\infty} |X|_k.$$

\begin{definition}
    The zeta functions of $X$ is
    $$Z(X,T) = \exp\left( \sum_{k=1}^{\infty} \frac{T^k}{k} \#X(\ff_{q^k}) \right)$$
    $$= \prod_{k=1}^{\infty} \frac{1}{(1-T^k)^{\#|X|_k}} \in 1+T\zz[[T]].$$
\end{definition}

\begin{question}
What does $Z(X,T)$ look like?
\end{question}
The answer was proposed by Andr\'e Weil in his celebrated Weil
conjectures. More precisely, Dwork \cite{Dw1} proved that $Z(X,T)$
is a rational function. Deligne \cite{D2} proved that the
reciprocal zeros and poles of $Z(X,T)$ are Weil $q-$numbers.

\subsection*{Moment Zeta Functions}
Let $f:X \mapsto Y /\ff_q$.  One has

$$X(\overline{\ff_q}) = \coprod_{y \in Y(\overline{\ff_q})} f^{-1}(y)(\overline{\ff_q}).$$
Similarly
$$X(\ff_q) = \coprod_{y \in Y(\ff_q)} f^{-1}(y)(\ff_q).$$
From this we get
$$\#X(\ff_{q^k}) = \sum_{y \in Y(\ff_{q^k})} \#f^{-1}(y)(\ff_{q^k})$$
for $k=1,2,3,\ldots$.  This number is known as the first moment of $f$ over $\ff_{q^k}$.

\begin{definition}
\label{dth_moment}
For $d \in \zz_{>0}$, the $d$-th moment of $f$ over $\ff_{q^k}$ is

$$M_d(f \otimes \ff_{q^k}) = \sum_{y \in Y(\ff_{q^k})} \#f^{-1}(y)(\ff_{q^{dk}})$$
$k=1,2,3,\ldots$
\end{definition}
\begin{definition}
The $d$-th moment zeta function of $f$ over $\ff_q$ is
$$
\begin{array}{lll}
Z_d(f,T) & = & \exp \left( \sum_{k=1}^{\infty} \frac{T^k}{k} M_d(f \otimes \ff_{q^k}) \right) \\
&=&\prod_{y \in |Y|} Z \left( f^{-1}(y)
\otimes_{\ff_{q^{\deg(y)}}} \ff_{q^{d\times \deg(y)}} ,
T^{\deg(y)} \right) \in 1+T\zz[[T]].

\end{array}
$$
\end{definition}

\begin{figure}[ht]
\caption{$f^{-1}(y)$}
\label{f_inv}
$$
\begin{xy} <1cm,0cm>:
(0.2,1.8)*{X};
(5.2,0)*{Y};
(0,0);(5,0);**@{-};
(2,0)="h1"*{\bullet};
"h1"+(0.4,0.2)*{y_1};
"h1"+(0,-2.6)+(0.2,0.2)*{f^{-1}(y_1)};
 "h1"+(0,-2);"h1"+(0,2);
**\crv{"h1"+(-0.5,0)&"h1"+(0.5,0)};
(3.5,0)="h2"*{\bullet};
"h2"+(0.4,0.2)*{y_2};
"h2"+(0,-2.6)+(0.2,0.2)*{f^{-1}(y_2)};
 "h2"+(0,-2);"h2"+(0,2);
**\crv{"h2"+(-0.5,0)&"h2"+(0.5,0)};
\end{xy}
$$
\end{figure}
Geometrically $M_d(f \otimes \ff_{q^k})$ can be thought of as certain point counting along the fibres of $f$.
Note that $M_d(f,k)$ will increase as $d$ increases.  Figure \ref{vary_d} illustrates this.
\begin{figure}[ht]
\caption{$f^{-1}(y)(\ff_{q^d})$}
\label{vary_d}
$$
\begin{xy} <0.75cm,0cm>:
(0.2,1.8)*{X};
(5.2,0)*{Y};
(2.6,-2.3)*{d=1};
(0,0);(5,0);**@{-};
(2,0)="h1"*{\bullet};
 "h1"+(0.15,-1);"h1"+(-0.15,1);
**\crv{"h1"+(-0.25,0)&"h1"+(0.25,0)};
(3.5,0)="h2"*{\bullet};
 "h2"+(0.15,-1);"h2"+(-0.15,1);
**\crv{"h2"+(-0.25,0)&"h2"+(0.25,0)};
"h1"+(-0.15,1);"h2"+(-0.15,1);**@{-};
"h1"+(0.15,-1);"h2"+(0.15,-1);**@{-};
\end{xy}
~~~
\begin{xy} <0.75cm,0cm>:
(0.2,1.8)*{X};
(5.2,0)*{Y};
(2.6,-2.3)*{d=2};
(0,0);(5,0);**@{-};
(2,0)="h1"*{\bullet};
 "h1"+(0,-2);"h1"+(0,2);
**\crv{"h1"+(-0.5,0)&"h1"+(0.5,0)};
(3.5,0)="h2"*{\bullet};
 "h2"+(0,-2);"h2"+(0,2);
**\crv{"h2"+(-0.5,0)&"h2"+(0.5,0)};
"h1"+(-0,2);"h2"+(0,2);**@{-};
"h1"+(0,-2);"h2"+(0,-2);**@{-};
"h1"+(-0.15,1);"h2"+(-0.15,1);**@{.};
"h1"+(0.15,-1);"h2"+(0.15,-1);**@{.};
\end{xy}
$$
As $d$ increases the area where we count points will also increase.
\end{figure}
The sequence of moment zeta functions $Z_d(f, T)$ measures the
arithmetic variation of rational points along the fibres of $f$.
It naturally arises from the study of Dwork's unit root conjecture
\cite{WAN05}.
\begin{question}~
    \begin{enumerate}
        \item For a given $f$, what is $Z_d(f,T)$?
        \item How does $Z_d(f,T)$ vary with $d$?
    \end{enumerate}
\end{question}

\subsection*{Partial Zeta Functions}

\begin{figure}[ht]
\caption{$f:X \mapsto X_1 \times \ldots \times X_n$}
$$
\xymatrix@R=3pt
{
f:X \ar[r]^{f_1} \ar@{.>}[dr]|{f_i} \ar[ddr]_{f_n} & X_1 \POS[];[d]**\dir{.}, \\
              & X_i  \POS[];[d]**\dir{.}\\
              & X_n
}
$$
\end{figure}
Assume $f:X \mapsto X_1 \times \ldots \times X_n$ defined by $x \mapsto (f_1(x), \ldots, f_n(x))$ is an embedding.
There are many ways to satisfy this property.  For example the addition of the identity function $f_n:X \mapsto X$ will assure $f$ is an embedding.

Let $d_1, \ldots, d_n \in \zz_{>0}$.  For $k=1,2,3,\ldots$, let

$$M_{d_1,\ldots, d_n} (f \otimes \ff_{q^k}) :=
\#\{x \in X(\overline{\ff_q}) | f_1(x) \in X_1(\ff_{q^{d_1 k}}),\ldots, f_n(x) \in X_n(\ff_{q^{d_n k}})\} < \infty$$

\begin{definition}
    Define the partial zeta function of $f$ over $\ff_q$ to be
$$Z_{d_1,\ldots,d_n}(f,T) = \exp \left( \sum_{k=1}^{\infty} \frac{T^k}{k} M_{d_1,\ldots,d_n}(f \otimes \ff_{q^k}) \right).$$
\end{definition}
The partial zeta function measures the distribution of rational
points of $X$ independently along the fibres of the $n$-tuple of
morphisms $(f_1, \cdots, f_n)$.

\begin{example}
    If $f_1:X \mapsto X_1$ and  $f_2= \Id:X \mapsto X$ then $Z_{1,d}(f,T) = Z_d(f_1,T)$.
\end{example}

Thus, partial zeta functions are generalizations of moment zeta functions.

\begin{question}~
    \begin{enumerate}
        \item What is $Z_{d_1,\ldots,d_n}(f,T)$?
        \item How does $Z_{d_1,\ldots,d_n}(f,T)$ vary as $\{d_1,\ldots, d_n\}$ varies?
    \end{enumerate}
\end{question}
We have
\begin{theorem}[\cite{WAN03}] The partial zeta function $Z_{d_1,\ldots,d_n}(f,T)$ is a rational function.
Furthermore, its reciprocal zeros and poles are Weil $q$-numbers.
\end{theorem}

\section{General properties of $Z(f,T)$.}

\subsection*{Trace Formula}  By Grothendieck \cite{Gr},
$Z(X,T)$ can be expressed in terms of $l$-adic cohomology. More
precisely, let $\overline{X} = X \otimes_{\ff_q}
\overline{\ff_q}$. Then,

\begin{theorem}
There are finite dimensional vector spaces $H_c^i(X)$ with
invertible linear action by $\Frob_q$ such that
$$Z(X,T) = \prod_{i=0}^{2 \dim(X)} \det(I-\Frob_q^{-1}T | H_c^i(X))^{(-1)^{i-1}},$$
where
$$H_c^i(X)= \begin{cases}
 H_c^i( \overline{X}, \qq_l ) & l \neq p, prime \\
H_{rig,c}(X, \qq_p), & l = p.
\end{cases}$$
\end{theorem}
This is used to show that $Z(X,T) \in \qq(T)$.
One should note:
\begin{enumerate}
    \item $Z(X,T)$ is independent of the choice of $l$.
    \item $\det(I- \Frob_q^{-1}T | H_c^i(X))$ may depend on the choice of $l$ due to possible
    cancelation. The conjectural independence on $\ell$ is still
    open in general.
\end{enumerate}

\subsection*{Reimann Hypothesis}
Fix an embedding of $\overline{\qq_l} \hookrightarrow \cc$.  Let $b_i = \dim  H_c^i(X)$.  Consider the factorization
$$\det(I-\Frob_q^{-1}T | H_c^i(X))= \prod_{j=1}^{b_i} (1 - \alpha_{ij}T), \alpha_{ij} \in \cc.$$
The $\alpha_{ij}$'s are Weil $q$-numbers, that is,
\begin{enumerate}
    \item The $\alpha_{ij}$'s are algebraic integers over $\qq$.
    \item For $\sigma \in \Gal({\overline{\qq} / \qq}), |\alpha_{ij}| = |\sigma(\alpha_{ij})| = \sqrt{q}^{\omega_{ij}}$ for some integer $\omega_{ij}$, called the weight of $\alpha_{ij}$ with $ 0\leq \omega_{ij} \leq i$,$\forall j=1,\ldots b_i$.
\end{enumerate}
The $l \neq p$ case was proved by Deligne \cite{D2} and the $l=p$
case by Kedlaya \cite{Ke}.

\subsection*{Slopes ($p$-adic Reimann Hypothesis)}

Consider an embedding $\overline{\qq_l} \hookrightarrow \cc_p$.  Then what is the $\ord_q(\alpha_{ij})\in \qq_{\geq 0}$?  This is referred to as the slope of $\alpha_{ij}$.

By Riemann Hypothesis,
$$\alpha_{ij} \overline{\alpha_{ij}} = q^{\omega_{ij}},$$

$$0 \leq \ord_q(\alpha_{ij}) \leq \ord_q(\alpha_{ij} \overline{\alpha_{ij}} ) = \omega_{ij} \leq i,$$

Further, Deligne's integrality theorem implies that

$$i-\dim(X) \leq \ord_q(\alpha_{ij}).$$

\begin{question} Given $X/ \ff_q$, the following questions arise:
    \begin{enumerate}
        \item What is $b_{i,l}: = b_{i}$?
        \item What is $\omega_{ij}$?
        \item What is the slope $\ord_q(\alpha_{ij})$?
    \end{enumerate}
\end{question}

\begin{example}
If $X$ is a smooth projective variety over $\ff_q$, then:
    \begin{enumerate}
        \item $H_c^i(X)$ is pure of weight $i$, i.e. $\omega_{ij} = i$ for $1 \leq j \leq b_i$.  Thus $b_{i,l}$ is independent of $l$.
        \item The $q$-adic Newton polygon (NP) of $\det(I-\Frob_q^{-1}T | H_c^i(X)) \in \zz[[T]]$
        lies above the Hodge polygon of $H_c^i(X)$.  This was conjectured by Katz \cite{Ka1} and proven by
        Mazur \cite{Ma} and Ogus \cite{BO}.  We will discuss this more later.
    \end{enumerate}
\end{example}


\section{Moment Zeta Functions}

Let $f:X \rightarrow Y/ \ff_q$.
For $d \in \zz_{>0}$, recall the $d$-th moment of $f \otimes \ff_{q^k}$ is
$$M_d(f \otimes \ff_{q^k}) = \sum_{y \in Y(\ff_{q^k})} \#f^{-1}(y)(\ff_{q^{dk}}).$$

\begin{question}~
    \begin{enumerate}
        \item How does $M_d(f \otimes \ff_{q^k})$ vary as $k$ varies?
        \item How does $M_d(f \otimes \ff_{q^k})$ vary with $d$?
        \item How does $M_d(f \otimes \ff_{q^k})$ vary with both $d$ and $k$?
    \end{enumerate}
\end{question}

\begin{definition}
    Define the $d$-th moment zeta function of $f$ to be
    $$Z_d(f,T) = \exp \left( \sum_{k=1}^{\infty} \frac{T^k}{k} M_d(f \otimes \ff_{q^{k}})\right).$$
\end{definition}
Observe for $d=1$ we have $Z_1(f,T) = Z(X,T)$.
Recall that $Z_d(f,T) \in \qq(T)$ and its reciprocal zeros and poles are Weil $q$-numbers.
This follows from the following more precise cohomological formula.

\begin{theorem} Let $l \neq p$.
    Let $\frF^i = R^i f_!\qq_l$ be the $i$-th relative $l$-adic cohomology with compact support.
Then $Z_d(f,T) =$
    $$\prod_{i=0}^{2 \dim(X/Y)} \prod_{j=0}^{d} \prod_{k=0}^{2 \dim(Y)} \det \left(I - \Frob_q^{-1} T | H_c^k( \overline{Y} , Sym^{d-j} \frF^i \otimes \bigwedge^j \frF^i) \right) ^{(-1)^{i+j+k-1}(j-1)}$$
\end{theorem}

\begin{proof} For an $l$-adic sheaf $\frF$ on $Y$, let $L(\frF, T)$ denote the L-function of $\frF$.
The trace formula in \cite{Gr} applies to the L-function $L(\frF,
T)$:
$$L(\frF,T) = \prod_{i=0}^{2 \dim(Y)} \det(I-\Frob_q^{-1}T | H_c^i(\bar{Y}, \frF))^{(-1)^{i-1}}.$$
The $d$-th Adams operation of a sheaf $\frF$ can be written as the
virtual sheaf \cite{WAN99}
$$[\frF]^d = \sum_{j \geq 0} (-1)^j (j-1) \left[Sym^{d-j} \frF \otimes \bigwedge^j \frF \right].$$
It follows that
$$\begin{array}{lll}
Z_d(f,T) &= &\prod_{y \in |Y|} Z \left( f^{-1}(y) \otimes_{\ff_{q^{\deg(y)}}} \ff_{q^{\deg(y)d}},T^{\deg(y)} \right) \\
    & = &\prod_{y \in |Y|} \prod_{i \geq 0} \det\left( I - T (\Frob_{q^{\deg(y)}}^{-1})^d T^{\deg(y)} | \frF_y^i \right)^{(-1)^{i-1}} \\
    & = &\prod_{i \geq 0} \prod_{y \in |Y|} \det\left( I - T^{\deg(y)} (\Frob_{q^{\deg(y)}}^{-1}) | [\frF_y^i]^d \right)^{(-1)^{i-1}} \\
    & = &\prod_{i \geq 0} L([\frF^i]^d / Y,T)^{(-1)^{i}}\\
    & = &\prod_{i \geq 0} \prod_{j \geq 0} L\left(\Sym^{d-j} \frF^i \otimes \bigwedge^j \frF^i, T\right)^{(-1)^{i+j}(j-1)}
\end{array}$$
  $ = \prod_k \prod_{i \geq 0} \prod_{j \geq 0} \det\left(I-T \Frob_q^{-1} | H_c^k(\overline{Y}, \Sym^{d-j} \frF^i \otimes \bigwedge^j \frF^i,T)\right)^{(-1)^{i+j+k-1}(j-1)}.$
\end{proof}
To use this formula one needs to know:
    \begin{enumerate}
     \item The total degree of $Z_d(f,T)$: number of zeros + number of
     poles.
     \item The high weight trivial factor which gives the main term.
     \item The vanishing of nontrivial high weight term which gives a good error bound.
    \end{enumerate}
Note:
\begin{enumerate}
 \item  There is an explicit upper bound for the total degree of $Z_d(f,T)$, which grows exponentially in
 $d$, see \cite{FW04}.
 \item There exists a total degree bound of the form $c_1 d^{c_2}$ which is a polynomial in $d$.
However, the constant $c_1$ is not yet known to be effective if $
\dim Y \geq 2$, see \cite{FW04}.
\end{enumerate}
\begin{question}
  How do we make $c_1$ effective?
\end{question}

\subsection*{Example: Artin-Schreier hypersurfaces}

Let
$$g(x_1,\ldots, x_n,y_1,\ldots,y_{n'}) \in \ff_q[x_1,\ldots, x_n,y_1,\ldots, y_{n'}].$$  We may also rewrite
this as $g= g_m + g_{m-1} +\ldots+ g_0$, where $g_i$ is the homogeneous part of degree $i$ and $g_m \neq 0$.

Consider:
$$\begin{array}{ll}
X: & \{x_0^p - x_0 =  g(x_1,\ldots, x_n,y_1,\ldots,y_{n'})\} \hookrightarrow  \aa^{n+n'+1} \\
Y: & \aa^{n'} \\
f: & X \mapsto Y, (x_0,x_1,\ldots,x_n, y_1,\ldots, y_{n'}) \mapsto (y_1,\ldots,y_{n'})
\end{array}
$$
One may then ask:
$$ M_d(f) = \#\{x_i \in \ff_{q^d},y_i \in \ff_q | x_0^p - x_0 = g(x,y)\}=?$$

Ideally for nice $g$, one hopes:
$$M_d(f) = q^{dn+n'} + O(q^{(dn+n')/2})$$

\begin{theorem}[Deligne, \cite{D1}]
Assume that $g$ is a Deligne polynomial of degree $m$, i.e., the
leading form $g_m$ is a smooth projective hypersurface in
$\pp^{n+n'}$ and $ p \nmid m$.  Then
    $$|M_1(f)-q^{n+n'}| \leq (p-1)(m-1)^{n+n'} q^{\frac{n+n'}{2}}.$$
\end{theorem}
For $d>1$, a similar estimate can be obtained in some cases.

Assume $f^{-1}(y)$ is a Deligne polynomial of degree $m$ for all
$y \in \aa^{n'}(\ff_{q^d})$.  Then, applying Deligne's estimate
fibre by fibre, one deduces
$$ \# f^{-1}(y)(\ff_{q^d}) = q^{dn} + E_y(d),$$
$$ |E_y(d)| \leq (p-1)(m-1)^n q^{\frac{dn}{2}},$$
where $E_y(d)$ is some error term.  From this, we get
$$
\begin{array}{lll}
M_d(f) & = & \sum_{y \in \aa^{n'}(\ff_q)} \#f^{-1}(y)(\ff_{q^d}) \\
& = & q^{dn+n'} + \sum_{y \in \aa^{n'}(\ff_q)} E_y(d).
\end{array}$$

Thus,  we get the ``trivial" estimate:
$$|M_d(f)-q^{dn+n'}| \leq (p-1)(m-1)^n q^{\frac{dn}{2}+n'}.$$
Ideally, one would hope to replace $n'$ by $n'/2$ in the above error bound.

If one applies the Katz type estimate via monodromy calculation as
in \cite{Ka2}, one gets $\sqrt{q}$ savings in good cases,  i.e.,
with error term $O(q^{\frac{dn}{2}+n' - \frac{1}{2}})$. This is
still far from the expected error bound $O(q^{\frac{dn+n'}{2}})$
if $n'\geq 2$.
\begin{definition}
The $d$-th fibered sum of $g$ is
    $$\bigoplus_Y^d g = g(x_{1 1}, \ldots , x_{1 n}, y_{1}, \ldots , y_{n'}) + \ldots + g(x_{d 1}, \ldots , x_{d n}, y_{1}, \ldots , y_{n'})$$
Observe the $y_i$ values remain the same while the $x_{i j}$ values vary.
\end{definition}

\begin{theorem}[Fu-Wan, \cite{FW04}]
 Assume $\bigoplus_Y^d g$ is a Deligne polynomial of degree $m$.  Then
\begin{enumerate}
 \item $|M_d(f) - q^{dn+n'}| \leq (p-1)(m-1)^{dn+n'} q^{\frac{dn+n'}{2}}$
 \item $|M_d(f) - q^{dn+n'}| \leq c(p,n,n')d^{3(m+1)^n-1} q^{\frac{dn+n'}{2}}$
\end{enumerate}
The constant $c$ is not known to be effective if $n' \geq 2$.
\end{theorem}

If $p$ does not divide $d$, then $\bigoplus_Y^d g$ is a Deligne
polynomial for a generic $g$ of degree $m$. Thus, the assumption
is satisfied for many $g$ if $p$ does not divide $d$. However, if
$p\mid d$, there are no such $g$.
\begin{question}
If $ p | d$, what would be the best estimate for $M_d(f)$?
\end{question}

\subsection{Example: Toric Calabi-Yau hypersurfaces}

This geometric example is studied in a joint work with A.
Rojas-Leon \cite{RW07}. Let $n \geq 2$. We consider
$$X: \{ x_1+ \ldots + x_n + \frac{1}{x_1\ldots x_n} -y = 0\} \hookrightarrow \GG_m^n\times \aa^1,$$
$$ Y = \aa^1,$$
$$f: (x_1,\cdots, x_n, y) \longrightarrow y.$$
  For $y \neq (n+1) \zeta$, with $\zeta^{n+1} = 1$ we have
  $$f^{-1}(y): x_1 + \ldots + x_n + \frac{1}{x_1 \ldots x_n} -y =0$$
  is an affine Calabi-Yau hypersurface in $\GG_m^n$.

For $n=2$ we have an elliptic curve.  For $n=3$ we have a K3 surface.  For $n=4$ we have a Calabi-Yau 3-fold.  Recall

$$M_d(f) = \sum_{y \in \ff_q} \#f^{-1}(y)(\ff_{q^d}).$$

For $d=1$, we have $M_1(f) = \#X(\ff_q) = (q-1)^n$. For every $y \in
\ff_q$, we have
$$\#f^{-1}(y)(\ff_{q^d}) = \frac{(q^d-1)^n-(-1)^n}{q^d}+ E_y(d),$$
where $E_y(d)$ is some error term with $|E_y(d)| \leq n q^{d(n-1)/2}$.  Thus,
$$M_d(f) = q \frac{(q^d-1)^n-(-1)^n}{q^d}+ \sum_{y \in \ff_q} E_y(d).$$
From this, we obtain the ``trivial" estimate
$$|M_d(f) -  \frac{(q^d-1)^n-(-1)^n}{q^{d-1}}| \leq n q^{d(n-1)/2+1}.$$

\begin{theorem}[Rojas-Leon and Wan, \cite{RW07}]
 If $p \nmid (n+1)$, then
\begin{enumerate}
 \item $$|M_d(f) -  \left( \frac{(q^d-1)^n-(-1)^n}{q^{d-1}}+ \frac{1}{2}(1+(-1)^d)q^{d(n-1)/2 +1}\right)
  \leq D q^{d(n-1)/2+\frac{1}{2}},$$ where
$D$ is an explicit constant depending only on $n$ and $d$.
 \item The purity decomposition of $Z_d(f,T)$ is determined.
\end{enumerate}
\end{theorem}

\begin{question}
 How do $M_d(f)$ and $Z_d(f,T)$ vary with $d$?
\end{question}

\section{Zeta functions of fibres}
We continue with the previous example.
\begin{example} For $y\in \ff_q$,
let
$$f^{-1}(y) = x_1+ \ldots + x_n + \frac{1}{x_1 \ldots x_n} -y = 0
\hookrightarrow \GG_m^n.$$ This is singular when $y \in \{ (n+1)
\zeta | \zeta^{n+1}=1 \}$. This forms the mirror family of
$$\{
x_0^{n+1} + \ldots + x_n^{n+1} - yx_0 \ldots x_n = 0\}.$$
\end{example}

Let $p \nmid (n+1), y \in \ff_q \setminus \{ (n+1) \zeta | \zeta^{n+1}=1 \}$.  Then
$$ Z(f^{-1}(y)/\ff_q,T) = Z\left( \left\{ \frac{(q^k -1)^n - (-1)^n}{q^k} \right\}_{k=1}^{\infty} ,T \right)
P_y(T)^{(-1)^n},$$ where $P_y(T) \in 1+T\zz[T]$ of degree $n$,
pure of weight $(n-1)$. Write
$$P_y(T) = (1 - \alpha_1(y)T) \ldots (1 - \alpha_n(y)T), \ \ |\alpha_i(y)| = \sqrt{q^{n-1}}.$$
Then we get the following:
\begin{cor}
 $$|\#f^{-1}(y)(\ff_q) - \frac{(q-1)^n- (-1)^n}{q}| \leq n\sqrt{q^{n-1}}.$$
\end{cor}

The star decomposition in \cite{WAN93}\cite{WAN04} implies

\begin{theorem}
There is a nonzero  polynomial $H_p(y) \in \ff_p[y]$ such that if
$H_p(y) \neq 0$ for some $y\in \ff_q$, then $\ord_q(\alpha_i(y)) =
i-1$ for $1 \leq i \leq n$.
\end{theorem}

Equivalently,  this family of polynomials $f^{-1}(y)$ is generically
ordinary. An alternative proof can be found in Yu \cite{Yu}.

\subsection*{Moment Zeta Functions}
For $d > 0$, recall
$$M_d(f) = \sum_{ y \in \ff_q} \#f^{-1}(y)(\ff_{q^d}),$$
$$M_d(f \otimes \ff_{q^k}) = \sum_{y \in \ff_{q^k}} \#f^{-1}(y)(\ff_{q^{dk}}), k=1,2,3,\ldots,$$
$$Z_d(f,T) = \exp\left( \sum_{k=1}^{\infty} \frac{T^k}{k} M_d(f \otimes \ff_{q^k}) \right) \in \qq(T).$$

Let
$$S_d(T)=
\prod_{k=0}^{[\frac{n-2}{2}]}\frac{1-q^{dk}T}{1-q^{dk+1}T}
\prod_{i=0}^{n-1}(1-q^{di+1}T)^{(-1)^{i+1}{n \choose i+1}}.$$

\begin{theorem} Assume that $(n+1)$ divides $(q-1)$. Then, the
$d$-th moment zeta function has the following factorization
$$Z_d({\mathbb A}^1, X_{\lambda})^{(-1)^{n-1}}
=P_d(T)({Q_d(T)\over P(d,T)})^{n+1}A_d(T)S_d(T).
$$
We now explain each of the above factors. First, $P_d(T)$ is the
non-trivial factor which has the form
$$P_d(T)=\prod_{a+b=d, 0\leq b\leq n}P_{a,b}(T)^{(-1)^{b-1}(b-1)},$$
and each $P_{a,b}(T)$ is a polynomial in $1+T{\Bbb Z}[T]$, pure of
weight $d(n-1)+1$, whose degree $r$ is given explicitly and which
satisfies the functional equation
$$
P_d(T)=\pm T^r q^{(d(n-1)+1)r/2} P_d(1/q^{d(n-1)+1}T).
$$
Second, $P(d,T)\in 1+T{\Bbb Z}[T]$ is the $d$-th Adams operation of
the ``non-trivial" factor in the zeta function of a singular fibre
$X_{t}$, where $t=(n+1)\zeta_{n+1}$ and $\zeta_{n+1}=1$. It is a
polynomial of degree $(n-1)$ whose weights are completely
determined. Third, the quasi-trivial factor $Q_d(T)$ coming from a
finite singularity has the form
$$Q_d(T)=\prod_{a+b=d, 0\leq b\leq n}Q_{a,b}(T)^{(-1)^{b-1}(b-1)},$$
where $Q_{a,b}(T)$ is a polynomial whose degree $D_{n,a,b}$ and
the weights of its roots are given. Finally, the trivial factor
$A_d(T)$ is given by
$$A_d(T)= (1-q^{\frac{d(n-1)}{2}}T)(1-q^{\frac{d(n-1)}{2}+1}T)
(1-q^{\frac{d(n-2)}{2}+1}T)$$ if $n$ and $d$ are even,
$$A_d(T)=(1-q^{\frac{d(n-2)}{2}+1}T)$$
if $n$ is even and $d$ is odd,
$$A_d(T)=(1-q^{\frac{d(n-1)}{2}}T)$$
if $n$ and $d$ are odd,
$$A_d(T)= (1-q^{\frac{d(n-1)}{2}+1}T)^{-1}$$
if $n$ is odd and $d$ is even.
\end{theorem}

\begin{cor}
Let $n=2$ and $f: \{ x_1+ x_2 + \frac{1}{x_1 x_2} -y=0\} \mapsto
y$ with $p \nmid 3$. Then,
$$Z_d(f,T)^{-1} = A_d(T) \frac{R_d(T)}{R_{d-2}(qT)},$$
where $A_d(T)$ is a trivial factor and $R_d(T) \in 1+T\zz[T]$ is a
non-trivial factor which is pure of weight $d+1$ and degree
$2(d-1)$.
\end{cor}
For all $d \leq 1$, $R_d(T)=1$.  $R_2(T)$ is a polynomial of degree $2$ and weight $3$.
This suggests that $R_2(T)$ comes from a rigid Calabi-Yau variety.  In general $R_d(T)$ is of weight $d+1$ and degree $2(d-1)$.

As always, we may ask what are the slopes of $R_d(T)$?

The above one parameter family of Calabi-Yau hypersurfaces is the
only higher dimensional example for which the moment zeta functions
are determined so far. It shows that the calculation of the moment
zeta function can be quite complicated in general. A related example
is the one parameter family of higher dimensional Kloosterman sums,
see \cite{FW05}\cite{FW07} for the L-function of higher symmetric
power which gives the main piece of the moment zeta function.

\subsection{$l$-adic Moment Zeta Function ($l \neq p$)}

Fix a prime $l \neq p$.  Given $\alpha \in \zz_l^*$ and $d_1 \equiv d_2 ~ \mod ~ (l-1)l^{k-1}$ for some $k$, then $\alpha^{d_1} \equiv \alpha^{d_2} ~ \mod ~ l^k$.

By rationality of $Z(f^{-1}(y),T)$ it follows that
$$\#f^{-1}(y)(\ff_{q^d}) = \sum_i \alpha_i(y)^d - \sum_j \beta_j (y)^d$$
for some $l$-adic integers $\alpha_i(y)$ and $\beta_j(y)$.
Consider
$$M_d(f) = \sum_{y \in Y(\ff_q)} \#f^{-1}(y)(\ff_{q^d}).$$

This can be rewritten as
$$ = \sum_{y \in Y(\ff_q)} \left( \sum_i \alpha_i(y)^d - \sum_j \beta_j (y)^d \right).$$
We may take some $D_l(f) \in \zz_{>0}$ such that if $d_1 \equiv d_2 ~ \mod ~ D_l (f) l^{k-1}$ then
\begin{enumerate}
    \item $M_{d_1} (f) \equiv M_{d_2}(f) ~ \mod~ l^k$.
    \item $Z_{d_1}(f,T) \equiv Z_{d_2}(f,T) ~\mod ~ l^k \in 1+ T\zz[[T]]$.
\end{enumerate}

\begin{definition} The $l$-adic weight space is defined to be
    $$W_l(f) = \left( \zz / D_l(f)\zz \right) \times \zz_l.$$

Let $s = (s_1,s_2) \in W_l(f)$.  Take a sequence of $d_i \in \zz_{>0}$ such that

\begin{enumerate}
    \item $d_i \rightarrow \infty$ in $\cc$,
    \item $d_i \equiv s_1 ~\mod~D_l(f)$,
    \item $d_i \rightarrow s_2 \in \zz_l$.
\end{enumerate}

With this we may define the $l$-adic moment zeta function
$$\zeta_s(f,T) = \lim_{i \rightarrow \infty} Z_{d_i}(f,T) \in 1 + T \zz_l[[T]].$$
This function is analytic in the $l$-adic open unit disk $|T|_l < 1$.
\end{definition}

\begin{question}
Is $\zeta_s(f,T)$ analytic in $|T|_l \leq 1$?  What about in $|T|_l < \infty$?
\end{question}

Embed $\zz$ in $W_l(f)$ in the following way:
$$\zz \hookrightarrow W_l(f),$$
$$ d \mapsto (d,d).$$

\begin{prop}
If $d \in \zz_{>0} \hookrightarrow W_l(f)$ then $\zeta_d(f,T) = Z_d(f,T) \in \qq(T)$.
\end{prop}

\begin{question}
What if $s \in W_l (f) \setminus \zz$?  This is open even when $f$ is a non-trivial family of elliptic curves over $\ff_p$.
\end{question}

\subsection{$p$-adic moment zeta functions ($l=p$)}

As in the $l$-adic case, one has a $p$-adic continuous result.

If $d_1 \equiv d_2 ~\mod~D_p(f) p^{k-1}, d_1 \geq d_2 \geq c_f k$
for some $k$ and sufficiently large constant $c_f$ then
$$M_{d_1}(f) \equiv M_{d_2}(f) ~ \mod ~ p^k.$$
Also, define in the same way as above
$$\zeta_{s,p}(f,T) = \lim_{i \rightarrow \infty} Z_{d_i}(f,T) \in 1 + T\zz_p[[T]].$$
As before consider the embedding:
$$\zz \hookrightarrow W_p(f),$$
$$ d \mapsto (d,d).$$

The following result was conjectured by Dwork \cite{Dw2}.

\begin{theorem}[Wan, \cite{WAN99}\cite{WAN1}\cite{WAN2}]
If $s=d \in \zz \hookrightarrow W_p(f)$, then $\zeta_{d,p}(f,T)$ is $p$-adic meromorphic in $|T|_p < \infty$.\
\end{theorem}
Furthermore,  we have
\begin{theorem}[\cite{WAN2}]
    Assume the $p$-rank $\leq 1$.   Then for each $s \in W_p(f)$, $\zeta_{s,p}(f,T)$ is $p$-adic meromorphic in $|T|_p < \infty$.
\end{theorem}
This can be extended a little further as suggested by Coleman.

\begin{theorem}[Grosse-Kl\"onne, \cite{GK}] Assume the $p$-rank $\leq 1$.
For $s=(s_1,s_2)$ with $s_1 \in \zz/ D_p(f)$ and $s_2 \in \zz_p /
p^{\epsilon}$ (small denominator), then $\zeta_{s,p}(f,T)$ is
$p$-adic meromorphic in $|T|_p < \infty$.
\end{theorem}

\begin{question}
In the case $s\in W_p(f)-\zz$ and $p$-rank $>1$, it is unknown if
$\zeta_{s,p}(f,T)$ is $p$-adic meromorphic, even on the closed unit
disk $|T|_p\leq 1$.
\end{question}

\section{Moment Zeta functions over $\zz$}
Consider
$$f:X \mapsto Y / \zz[\frac{1}{N}].$$
The $d$-th moment zeta function of $f$ is:
    $$\zeta_d(f,s) = \prod_{p \nmid N} Z_d(f \otimes \ff_p, p^{-s}).$$
Is this $\cc$-meromorphic in $s \in \cc$?
Is $\zeta_d(f,s)$ or its special values $p$-adic continuous in some sense?
If so, its $p$-adic limit $\zeta_s(f) (s \in \zz_p)$ is a $p$-adic zeta function of $f$.

\begin{example} Consider the map
$$f:\{ x_1+x_2 + \frac{1}{x_1 x_2} -y=0\} \mapsto y.$$

Then
$$Z_d(f \otimes \ff_p,T)^{-1} = A_d(T) \frac{R_d(f \otimes \ff_p,T)}{R_{d-2}(f \otimes \ff_p,pT)}$$
where $A_d(T)$ is a trivial factor and $R_d$ is a non-trivial
factor of degree $2(d-1)$ and weight $d+1$.

$$R_d(T) \leftrightarrow f^{\otimes d} = \{x_{11} + x_{12} + \frac{1}{x_{11} x_{12}} = \ldots = x_{d1} + x_{d2} + \frac{1}{x_{d1} x_{d2}}\}$$
\end{example}

\begin{example}
For $d=2$, we have
$$x_1+x_2+\frac{1}{x_1 x_2} = y_1+y_2+\frac{1}{y_1 y_2}.$$

As Matthias Schuett observed during the workshop, $R_2(T)
\leftrightarrow$ the unique new form of weight $4$ and level $9$.
\end{example}

\begin{conjecture}
Is $\prod_p R_d\left(f \otimes \ff_p, p^{-s}\right)$ meromorphic in $s \in \cc$ for all $d$.
\end{conjecture}

The conjecture is known to be true if $d\leq 2$. It should be
realistic to prove it for all positive integers $d$.

\section{$l$-adic Partial Zeta Functions}
We now consider the system of maps where $X \mapsto X_1 \times
\ldots \times X_n$ is an embedding (See Figure \ref{map_system}).

\begin{figure}[ht]
\caption{$f:X \mapsto X_1 \times \ldots \times X_n$}
\label{map_system}
$$
\xymatrix@R=3pt
{
f:X \ar[r]^{f_1} \ar@{.>}[dr]|{f_i} \ar[ddr]_{f_n} & X_1 \POS[];[d]**\dir{.}, \\
              & X_i  \POS[];[d]**\dir{.}\\
              & X_n
}
$$
\end{figure}

  This allows us to define the partial zeta function

$$Z_{d_1,\ldots,d_n}(f,T) = \exp\left( \sum_{k=1}^{\infty} \frac{T^k}{k} \#\{ x \in X(\overline{\ff_q}) | f_i(x) \in X_i(\ff_{q^{d_i k}}) \right) \in \qq(T).$$

\begin{question}
Is there any $p$-adic or $l$-adic continuity result as $\{d_1, \ldots ,d_n\}$ varies $p$-adically or $l$-adically?
\end{question}

\begin{example}
Consider the surface and three projection maps:
$$f:x_1+x_2+\frac{1}{x_2 x_2} -x_3=
\xymatrix@R=0pt
{
0 \ar[r]^{f_1} \ar[dr]|{f_2} \ar[ddr]_{f_3} & x_1\\
              & x_2\\
              & x_3
}$$
Thus
$$M_{d_1,d_2,d_3}(f) = \#\{(x_1,x_2,x_3) | x_1+ x_2 + \frac{1}{x_1 x_2} -x_3 =0, x_i \in \ff_{q^{d_i}},i=1,2,3\}$$
Is there a continuity result as $\{d_1, d_2, d_3 \}$ vary?
\end{example}

\section{Zeta Functions of Toric Affine Hypersurfaces}

Let $\triangle \subset \rr^n$ be an $n$-dimensional integral polytope.  Let $f \in \ff_q[x_1^{\pm1} ,\ldots , x_n^{\pm1}]$  with
$$f = \sum_{u \in \triangle \cap \zz^n} a_u X^u, a_u \in \ff_q$$
such that $\triangle(f) = \triangle$.  That is, $a_u \neq 0$ for
each $u$ which is a vertex of $\triangle$.

\begin{question}
Consider the toric affine hypersurface
$$U_f:\{f(x_1,\ldots,x_n)=0\} \hookrightarrow \GG_m^n.$$
    \begin{enumerate}
        \item $\#U_f(\ff_q) = $?
        \item $Z(U_f,T) = $?
    \end{enumerate}
\end{question}

\begin{definition}~
    \begin{enumerate}
        \item If $\triangle' \subset \triangle$ is a face of $\triangle$, define
        $$f^{\triangle'} = \sum_{u \in \triangle' \cap \zz^n} a_u X^u.$$
        \item $f$ is $\triangle$-regular if for every face $\triangle'$ (of any dimension) of $\triangle$, the system
        $$f^{\triangle'} = x_1 \frac{\partial f^{\triangle'}}{\partial x_1} = \ldots = x_n \frac{\partial f^{\triangle'}}{\partial x_n} = 0$$
        has no common zeros in $\GG_m^n(\overline{\ff_q})$.
    \end{enumerate}
\end{definition}

\begin{theorem}[GKZ, \cite{GKZ}]~
    \begin{enumerate}
        \item There is a nonzero polynomial $\disc_{\triangle} \in \zz[a_u | u \in \triangle \cap \zz^n]$ such that
        $f$ is $\triangle$-regular if and only if $\disc_{\triangle}(f)
\neq 0$ in $\ff_q$.  In other words, $\disc_{\triangle}$ is an
integer coefficient polynomial that will determine
$\triangle$-regularity.
        \item $\triangle(\disc_{\triangle})$ is determined.  This is referred to as the secondary polytope.
    \end{enumerate}
\end{theorem}

\begin{question}
    For which $p$,  $\disc_{\triangle} \otimes \ff_p \neq 0$?
\end{question}

\begin{definition}
Let $C(\triangle)$ be the cone in $\rr^{n+1}$ generated by $0$ and
$(1,\triangle)$.

\begin{figure}[ht]
\caption{$C(\triangle)$}
$$
\begin{xy} <1cm,0cm>:
(0,0)*{\bullet}; (5,0)**@{--};
(0,0); (5,2.5)**@{--};
(0,0); (4.5,-1)**@{--};
(0,-0.25)*{\text{0}};
(4,0)*{\bullet};
(4,2)*{\bullet}**@{-};
(3.6,-0.8)*{\bullet}**@{-};
(4,0)**@{-};
(3.6,-1.2)*{\text{(4,$4\triangle$)}};
(3,0)*{\bullet};
(3,1.5)*{\bullet}**@{-};
(2.7,-0.6)*{\bullet}**@{-};
(3,0)**@{-};
(2.7,-1)*{\text{(3,$3\triangle$)}};
(2,0)*{\bullet};
(2,1)*{\bullet}**@{-};
(1.8,-0.4)*{\bullet}**@{-};
(2,0)**@{-};
(1.8,-0.8)*{\text{(2,$2\triangle$)}};
(1,0)*{\bullet};
(1,0.5)*{\bullet}**@{-};
(0.9,-0.2)*{\bullet}**@{-};
(1,0)**@{-};
(0.9,-0.6)*{\text{(1,$\triangle$)}}
\end{xy}
$$
\end{figure}

\begin{enumerate}
\item Define
$$W_{\triangle}(k) = \#\{(k, k\triangle ) \cap \zz^{n+1}\},k=0,1,\ldots$$
The Hodge numbers of $\Delta$ are defined by
$$h_{\triangle}(k) = W_{\triangle}(k) - {n+1 \choose 1} W_{\triangle}(k-1) + {n+2 \choose 2} W_{\triangle}(k-2) - \ldots,$$
$$h_{\triangle}(k) = 0, \text{ if } k \geq n+1.$$
\item $\deg(\triangle) = d(\triangle) = n!Vol(\triangle) = \sum_{k=0}^{n} h_{\triangle}(k)$.
\end{enumerate}
\end{definition}

\begin{theorem}[Adolphson-Sperber \cite{AS}, Denef-Loesser \cite{DL}] Assume $f/\ff_q$ is $\triangle$-regular.  Then
    \begin{enumerate}
        \item $Z(U_f,T) = \prod_{i=0}^{n-1} (1 - q^i T)^{(-1)^{n-i} { n \choose i+1}} P_f(T)^{(-1)^n}$ with $P_f(T) \in 1 + T\zz[T]$ is of degree $d(\triangle)-1$.
        \item $P_f(T) = \prod_{i=1}^{d(\triangle)-1} ( 1 - \alpha_i(f) T), |\alpha_i(f)| \leq \sqrt{q}^{n-1}$. In particular,
        $$|\# U_f(\ff_q) - \frac{(q-1)^n - (-1)^n}{q}| \leq (d(\triangle)-1) \sqrt{q}^{n-1}.$$
        The precise weights of the $\alpha_i(f)$'s were also
        determined by Denef-Loesser.
    \end{enumerate}
\end{theorem}

\begin{question}
    For $i=1,2,\ldots, d(\triangle)-1$,  what is $\ord_q(\alpha_i(f))=$?
\end{question}

\section{Newton and Hodge Polygons}
Write
$$P_f(T) = 1 + c_1 T + c_2 T^2 +\ldots.$$
The $q$-adic Newton polygon of $P_f(T)$ is the lower convex closure in $\rr^2$ of the points $(k, \ord_q(c_k)),(k=0,1,\ldots, d(\triangle)-1)$. Denote this Newton polygon by $NP(f)$. Note that $NP(f)=NP(f\otimes \ff_{q^k})$ for
all $k$.

\begin{figure}[ht]
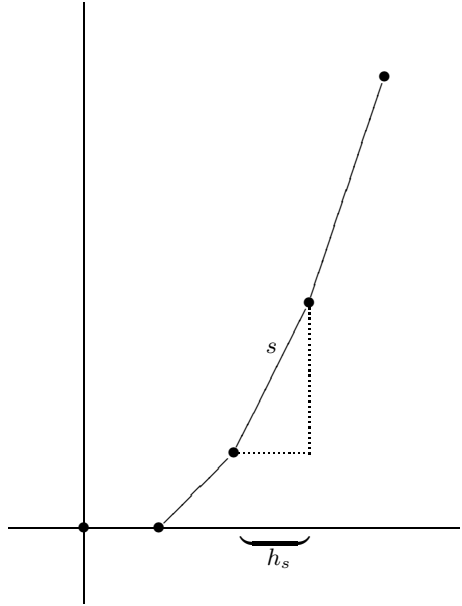

\caption{Newton Polygon}
$$
\xy<1cm,0cm>:
(0,-0.2)="o";
(1,-0.2)="h0";
(2,-0.2)="h1";
(3,-0.2)="h2";
(4,-0.2)="h3";
(0,0)*{\bullet};
(1,0)*{\bullet} **@{.};
(2,1)*{\bullet}="t1" **@{-};
(3,3)*{\bullet}="t2" **@{-};
(4,6)*{\bullet}="t3" **@{-};
(2,1);
(3,1)**@{.};
(3,3)**@{.};
"h1"+(0.4,0);"h2"**@{}**\frm{_)};
"h1"+(0.6,-0.2)*\txt\small{$h_s$};
(2.5,2.4)*\txt\small{$s$};
(0,-1);(0,7) **@{-};
(-1,0);(5,0) **@{-};
\endxy
$$
\end{figure}

\begin{prop} Let $h_s$ denote the horizontal length of the slope
$s$ side in $NP(f)$. Then, $P_f(T)$ has exactly $h_s$ reciprocal
zeros $\alpha_i(f)$ such that $ord_q(\alpha_i(f)) = s$ for each $s
\in \qq_{\geq 0}$.
\end{prop}

\begin{definition}
    The Hodge polygon of $\triangle$, denoted by $HP(\Delta)$,
     is the polygon in $\rr^2$ with a side of slope $k-1$ with horizontal length $h_{\triangle}(k)$ for $1\leq k\leq n$
     and vertices
    $$(0,0), \left( \sum_{m=1}^{k} h_{\triangle}(m), \sum_{m=1}^{k} (m-1)h_{\triangle}(m) \right), k=1,2,\ldots,n.$$
\end{definition}

\begin{theorem}[Adolphson-Sperber \cite{AS}]
    The $q$-adic Newton polygon lies above the Hodge polygon, i,e. $NP(f) \geq HP(\triangle)$.  In addition, the endpoints of the two coincide.
\end{theorem}

\begin{definition}
If $NP(f)=HP(\Delta)$, then $f$ is called ordinary.
\end{definition}

\begin{question}
    When is $f$ ordinary?  One hopes this is often.
\end{question}
Let
$$M_p(\triangle) = \{f\in \overline{\ff_p}[x_1^{\pm1}, \cdots, x_n^{\pm1}] |
\triangle(f) = \triangle, \text{$f$ is $\triangle$-regular} \}.$$

\begin{theorem}[Grothendieck, \cite{Ka2}]
There exists a generic Newton polygon, denoted by $GNP(\Delta, p)$,  such that
$$GNP(\triangle,p) = \inf\{NP(f) | f \in M_p(\triangle)\}$$
Hence for any $f \in M_p(\triangle)$,
$$NP(f) \geq GNP(\triangle,p) \geq HP(\triangle),$$
where the first inequality is an equality for most $f$ (generic
$f$).
\end{theorem}

\begin{question}
    Given $\triangle$, for which $p$,  is $GNP(\triangle,p) = HP(\triangle)$?  In other words, when is $f$ generically ordinary?
\end{question}
This suggests the following conjecture.
\begin{conjecture}[Adolphson-Sperber \cite{AS}]
    For each $p \gg 0$, $GNP(\triangle,p) = HP(\Delta)$.
\end{conjecture}
This is false in general.  Some counterexamples can be found in \cite{WAN93}.

\begin{definition}~
    \begin{enumerate}
        \item
        $$
        \begin{array}{lll}
        S(\triangle) &= &\text{ the semigroup }C(\triangle) \cap \zz^{n+1}. \\
        S_1(\triangle) &=&\text{ the semigroup generated by }(1,\triangle) \cap \zz^{n+1}.
        \end{array}
        $$
        \item Define the exponents of $\triangle$ as
        $$
        \begin{array}{lll}
        I(\triangle) &= \inf\{D>0 | Du \in S_1(\triangle) , \forall u \in S(\triangle)\} \\
        I_{\infty}(\triangle) &= \inf\{D>0 | Du \in S_1(\triangle) , \forall u \in S(\triangle), u \gg 0\}
        \end{array}
        $$
    \end{enumerate}
\end{definition}

\begin{conjecture}
If $p \equiv 1~ mod ~ I(\triangle)$ or if $p \equiv 1 ~mod~I_{\infty}(\triangle)$ for $p \gg0$,  then
\begin{enumerate}
\item $\disc_{\triangle} \otimes \ff_p \neq 0$,
 \item $GNP(\triangle,p) = HP(\triangle)$.
\end{enumerate}
\end{conjecture}
Part (2) is a weaker version of the conjecture in \cite{WAN93}.

\section{Generic Ordinarity}

\subsection*{Toric Hypersurface}

Let $\triangle\subset \rr^n$ be a $n$-dimensional integral polytope and $p$ a prime.  Let $d(\triangle) = n!Vol(\triangle)$.  Define
$$M_p(\triangle) = \{ f \in \overline{\ff_p}[x_1^{\pm1}, \ldots,  x_n^{\pm1}] | \triangle(f) = \triangle,
\text{$f$ is $\triangle$-regular}\}.$$

For each $ f \in M_p(\triangle)$, let $NP(f)$ be the Newton
polygon of the interesting factor $P_f(T)$ of the zeta function
$Z(U_f ,T)$.  Note that changing the ground field will not change
the Newton polygon. Recall that
$$NP(f) \geq GNP(\triangle,p) \geq HP(\triangle).$$

Note that $NP(f)$ is defined in a completely arithmetic fashion
and is dependent on the coefficients of the polynomial $f$. On the
other hand, $GNP(\triangle,p)$ is independent of coefficients
while $HP(\triangle)$ is obtained combinatorially. If
$GNP(\triangle,p) = HP(\triangle)$, we refer to $p$ as ordinary
for $\triangle$.
\begin{figure}[ht]
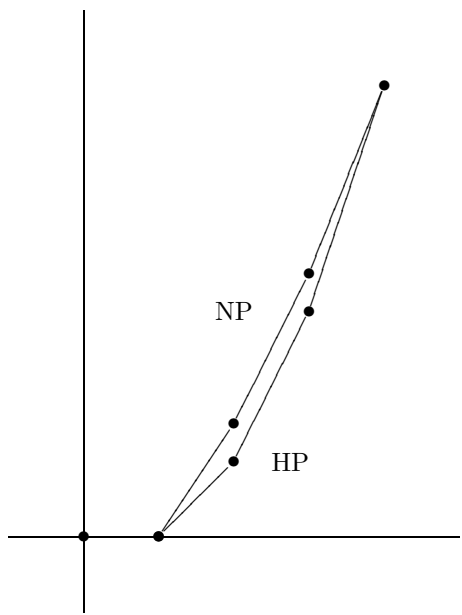

\caption{$NP \geq HP$}
$$
\xy<1cm,0cm>:
(0,-0.2)="o";
(1,-0.2)="h0";
(2,-0.2)="h1";
(3,-0.2)="h2";
(4,-0.2)="h3";
(0,0)*{\bullet};
(1,0)*{\bullet} **@{.};
(2,1)*{\bullet}="t1" **@{-};
(3,3)*{\bullet}="t2" **@{-};
(4,6)*{\bullet}="t3" **@{-};
(3,3.5)*{\bullet}="t3" **@{-};
(2,1.5)*{\bullet}="t3" **@{-};
(1,0)*{\bullet}="t3" **@{-};
(0,-1);(0,7) **@{-};
(-1,0);(5,0) **@{-};
(2,3)*{\txt{NP}};
(2.5,1)*!L{\txt{HP}};
\endxy
$$
\end{figure}
\begin{conjecture}[Adolphson-Sperber]
    For any $\triangle$, $p$ is ordinary for all $p \gg 0$.
\end{conjecture}

This conjecture is too strong as Example \ref{counterexample} illustrates.

\begin{example} 
\label{counterexample}
    Let $f = a_0+ a_1 x_1 + \ldots + a_n x_n + a_{n+1}x_1 x_2 \ldots x_n$ and
    $$\triangle = Conv((0,\ldots, 0),(1,\ldots,0),\ldots, (0, \ldots, 1), (1,1,\ldots,1)).$$  Therefore $d(\triangle) = n$ for $n \geq 2$. Furthermore, $\Delta$ is an empty simplex,
i.e., a simplex with no lattice points other than vertices. It follows that
    \begin{enumerate}
        \item $p$ is ordinary for $\triangle$ if and only if $p \equiv 1 ~\mod ~ (n-1)$. This implies
        \item If $n \geq 4$, then the Adolphson-Sperber conjecture is false.
    \end{enumerate}
\end{example}
\subsection*{Convex Triangulation}

\begin{definition}~
    \begin{enumerate}
        \item A triangulation of $\triangle$ is a decomposition
        $$\triangle = \bigcup_{i=1}^{m} \triangle_i,$$
        such that each $\triangle_i$ is a simplex, $\triangle_i \cap \triangle_j$ is a common face for both $\triangle_i$ and $\triangle_j$.
        \item The triangulation is called \bf convex \rm if there is a piecewise linear function $\phi:\triangle \mapsto \rr$ such that
        \begin{enumerate}
            \item $\phi$ is convex i.e. $\phi(\frac{1}{2} x+ \frac{1}{2} x') \leq \frac{1}{2} \phi(x) + \frac{1}{2} \phi(x')$, for all $x,x' \in \triangle$.
            \item The domains of linearity of $\phi$ are precisely the $n$-dimensional simplices $\triangle_i$ for $1 \leq i \leq m$.
        \end{enumerate}
    \end{enumerate}
\end{definition}

\begin{figure}[ht]
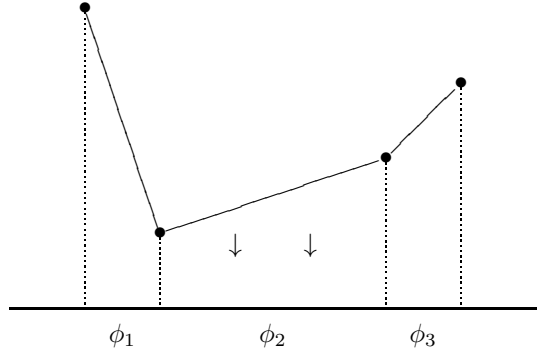

\caption{Piecewise projection down}
\label{proj_down}
$$\xy<1cm,0cm>:
(1,4)="v1";
(2,1)="v2";
(5,2)="v3";
(6,3)="v4";
"v1"*{\bullet};
"v2"*{\bullet} **@{-};
"v3"*{\bullet} **@{-};
"v4"*{\bullet}**@{-};
(0,0);(7,0) **@{-};
"v1"; (1,0) **@{.};
"v2"; (2,0) **@{.};
"v3"; (5,0) **@{.};
"v4"; (6,0) **@{.};
(1.5,-0.2)*!U{\phi_1};
(3.5,-0.2)*!U{\phi_2};
(5.5,-0.2)*!U{\phi_3};
(3,1)*!U{\downarrow};
(4,1)*!U{\downarrow};
\endxy$$
\end{figure}

\subsection*{Basic Decomposition Theorem}
There are several decompistion theorems for generic Newton polygon
in \cite{WAN93}\cite{WAN04}. They are based on a suitable convex
trangilation of $\Delta$ and thus they can be extended to prove the
following general decomposition theorem.

\begin{theorem}~
    \begin{enumerate}
        \item Let $\triangle= \cup_{i=1}^{m} \triangle_i$ be a convex integral triangulation of $\triangle$.
        If $p$ is ordinary for each $\triangle_i$, $1 \leq i \leq m$, then $p$ is ordinary for $\triangle$.

        \item If $\triangle$ is a simplex and $p \equiv 1~mod~ d(\triangle)$, then $p$ is ordinary.
    \end{enumerate}
\end{theorem}

\begin{cor}
    If $p \equiv 1 ~ mod ~( lcm(d(\triangle_1), \ldots , d(\triangle_m)))$, then $p$ is ordinary.
\end{cor}

\begin{example}
Let $A$ be the convex closure of (-1,-1), $(1,0)$ and $(0,1)$ in
$\rr^2$. The star decomposition in Figure \ref{A_simplex_other} is
convex and integral.
\begin{figure}[ht]
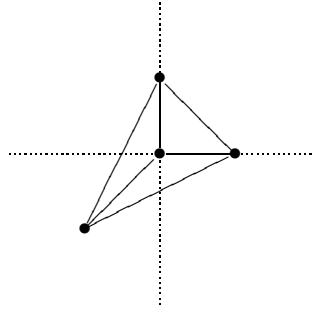

\caption{Star decomposition of $A$} \label{A_simplex_other}
$$\xy<1cm,0cm>:
(-2,0); (2,0); **@{.};
(0,-2); (0,2); **@{.};
(0,1)*{\bullet}; (1,0); **@{-};
(-1,-1)*{\bullet}; **@{-};
(-1,-1); (1,0)*{\bullet}; **@{-};
(0,0)*{\bullet};
(1,0); **@{-};
(0,1); **@{-};
(-1,-1); **@{-};
\endxy$$
\end{figure}
\end{example}
More generally,

\begin{example} Consider
    $$f:\{ x_1 + x_2 + \ldots + x_n + 1/x_1 x_2 \ldots x_n - y=0\} \text{ over } \ff_p.$$
    This is generically ordinary for all $p$. The proof uses the same star decomposition.
\end{example}

\begin{example}
    Let $\triangle = \{(d,0,\ldots,0),(0,d,0,\ldots,0), \ldots , (0,\ldots,d),(0,\ldots,0)\}$.
    We may make a parallel hyperplane cut as in Figure \ref{figparallel}.
    This will make $d(\triangle_i)=1$ for each piece $\triangle_i$ of the
    decomposition, see \cite{WAN93}.
    This proves that
the universal family of affine (or projective) hypersurfaces of
degree $d$ and $n$ variables over $\ff_p$ is also generically
ordinary for every $p$. The projective hypersurface (complete
intersection) case was first proved by Illusie \cite{Il}.

    \begin{figure}[ht]
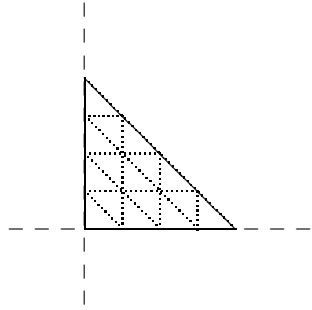


    \caption{Parallel Hyperplane Decomposition into simplices}
        \label{figparallel}
    $$
    \xy<1cm,0cm>:
(0.5,0); (0.5,1.5) **@{.};
(1,0); (1,1) **@{.};
(1.5,0); (1.5,0.5) **@{.};
(0,0.5); (1.5,0.5) **@{.};
(0,1); (1,1) **@{.};
(0,1.5); (0.5,1.5) **@{.};
(0,0.5); (0.5,0); **@{.};
(0,1); (1,0); **@{.};
(0,1.5); (1.5,0); **@{.};
(0,2); (2,0); **@{.};
(-1,0); (3,0) **@{--};
(0,-1); (0,3) **@{--};
(0,0) \PATH~={**@{-}} '(0,2)'(2,0)'(0,0);
    \endxy
$$
\end{figure}

\end{example}

\begin{cor}
    If $n = \dim(\triangle) = 2$, then $p$ is ordinary for $\triangle$ for all $p$.
\end{cor}

\begin{cor}
    If $n= \dim(\triangle) = 3$, then $p$ is ordinary for $p > 6 Vol(\triangle)$.
\end{cor}
This corollary is proven by showing stability of the $p$-action on the weight. This is a different argument than by proving $d(\triangle_i)=1$ argument.

\begin{definition}
Let $\triangle$ be an $n$-dimensional integral convex polytope in
$\rr^n$. Assume that $0$ (origin) is in the interior of $\triangle$.
Given such a situation, define $\triangle^* \subset \rr^n$ by:
$$\triangle^* = \{(x_1, \ldots, x_n) \in \rr^n | \sum_{i=1}^n x_iy_i \geq -1, ~\forall (y_1, \ldots ,y_n ) \in \triangle\}$$
\end{definition}
Observe $\triangle^*$ is also a convex polytope in $\rr^n$, though it may not have integral vertices.  Also observe $(\triangle^*)^*=\triangle$.

\begin{definition}
    $\triangle$ is called reflexive if $\triangle^*$ is also integral.
\end{definition}

\begin{cor}
    If $n= \dim(\triangle)=4$ and if $\triangle$ is reflexive, then $p$ is ordinary for $\triangle$
    for all $p>12 Vol(\Delta)$.
\end{cor}

\subsection{Slope Zeta Function}
The concept of slope zeta functions was developed for arithmetic
mirror symmetry as we will describe here. More information can be
found in \cite{WAN06}\cite{WAN07}.

Let $(X,Y)$ be a mirror pair over $\ff_q$.  Candelas, de la Ossa and
Rodriques-Villegas in \cite{COGL91} desired a possible mirror
relation of the type
$$Z(X,T) = \frac{1}{Z(Y,T)}$$
for $3$ dimensional Calabi-Yau varieties.  This is not true.  If
this were the case then
$$\sum \frac{T^k}{k} \#X(\ff_q) = \sum \frac{T^k}{k} (-\#Y(\ff_q)).$$
Therefore
$$ \#X(\ff_q) = - \#Y(\ff_q),$$
which is impossible for large $q$ on nonempty varieties.

The question is then to modify the zeta function suitably so that
the desired mirror relation holds.  The slope zeta function was
introduced for this purpose.

\begin{definition} Write $Z(X,T) = \prod_i (1 - \alpha_i T)^{\pm1} \in \cc_p(T)$.
    \begin{enumerate}
        \item The slope zeta function of $X$ is defined to be the
        following two variable function:
        $$S(X,U,T) = \prod_i (1 - U^{{\rm ord}_q(\alpha_i)}T)^{\pm1}.$$

        \item If $f: X \mapsto Y$ defined over $\ff_q$ (a nice family) then the slope zeta function of $f$ is the generic one among
        $S(f^{-1}(y),U,T)$ from all $y \in Y$, denoted by $S(f,U,T)$.
    \end{enumerate}

    \begin{conjecture}
        Let $X$ be a $3$-dimensional Calabi-Yau variety over $\qq$.  Assume that $X$ has a mirror over $\qq$.  Then the generic family containing $X$ as a member is generically ordinary for all $p \gg 0$.
    \end{conjecture}
    This conjecture implies the following

    \begin{conjecture}[Arithmetic Mirror Conjecture]
    Let $\{f,g\}$ be two generic mirror families of a $3$-dimensional Calabi-Yau variety over $\qq$.  Then
$$S(f \otimes \ff_p,U,T) = \frac{1}{S(g \otimes \ff_p,U,T)}$$
for all $p\gg 0$.
    \end{conjecture}
\end{definition}


\begin{thebibliography}{12}

\bibitem{AS}A. Adolphson and S. Sperber.  \it Exponential sums and Newton polyhedra: Cohomology and estimates\rm.
Ann. Math., 130 (1989), 367-406.

\bibitem{BO} P. Berthelot and A. Ogus,
\it Notes on Crystalline Cohomology, \rm Princeton University
Press, 1978.

\bibitem{COGL91}
P. Candelas, X. de la Ossa, F. Rodriques-Villegas, \it Calabi-Yau
manifolds over finitef fields II, \rm Fields Instit. of Commun.,
38(2003).


\bibitem{DL}J. Denef and F. Loesser, Weights of exponential sums, intersection cohomology, and Newton polyhedra,
Invent. Math., 106(1991), no.2, 275-294.

\bibitem{D1} P. Deligne,
\it Applications de la Formule des Traces aux Sommes
Trigonom\'etriques, \rm in Cohomologie \'Etale (SGA
$4\frac{1}{2}$), 168-232, Lecture Notes in Math. 569,
Springer-Verlag 1977.

\bibitem{D2} P. Deligne,
\it La Conjecture de Weil II, \rm Publ. Math. IHES 52 (1980),
137-252.



\bibitem{Dw1} B. Dwork,
\it On the rationality of the zeta function of an algebraic
variety, \rm Amer. J. Math., 82(1960), 631-648.

\bibitem{Dw2} B. Dwork,
\it Normalized period matrices II, \rm Ann. Math., 98(1973), 1-57.

\bibitem{FW04}L. Fu and D. Wan,
\it Moment L-functions, partial L-functions and partial exponential sums,
\rm Math. Ann., 328(2004), 193-228.

\bibitem{FW05} L. Fu and D. Wan,
\it L-functions for symmetric products of Kloosterman sums, \rm J.
Reine Angew. Math., 589(2005), 79-103.


\bibitem{FW07} L. Fu and D. Wan,
\it Trivial factors for L-functions of symmetric products of
Kloosterman sheaves, \rm Finite Fields \& Appl., to appear.

\bibitem{GKZ} I.M. Gelfand, M.M. Kapranov and A.V. Zelevinsky,
\it Discriminatns, Resultants and Multidimensional Determinants,
\rm Birkh\"user Boston, Inc., Boston, MA, 1994.

\bibitem{GK} E. Grosse-Kl\"onne,
\it On families of pure slope L-functions, \rm Documenta Math.,
8(2003), 1-42.

\bibitem{Gr} A. Grothendick,
\it Formule de Lefschetz et rationalit\'e des fonctions L, \rm
S\'eminaire Bourbaki, expos\'e 279, 1964/65.

\bibitem{Il}
L. Illusie. \it Ordinarit\'e des intersections compl\'etes
g\'en\'erales, \rm Grothendieck Festschrift, Vol. II (1990).
375-405.

\bibitem{Ka0} N. Katz,
\it On a theorem of Ax, \rm Amer. J. Math., 93(1971), 485-499.

\bibitem{Ka1} N. Katz,
\it Slope filtration of F-crystals, \rm Ast\'erisque 63(1979),
113-164.

\bibitem{Ka2} N. Katz,
\it Frobenius-Schur indicator and the ubiquity of Brock-Granville
quadratic excess, \rm Finite Fields \& Appl., 7(2001), 45-69.

\bibitem{Ke} K. Kedlaya,
\it Fourier transforms and p-adic ``Weil II", \rm Compositio
Mathematica, 142(2006), 1426-1450.

\bibitem{Ma} B. Mazur,
\it Frobenius and the Hodge filtration, \rm Bull. Amer. Math.
Soc., 78(1972), 653-667.

\bibitem{RW07}A. Rojas-Leon and D. Wan,
\it Moment zeta functions for toric Calabi-Yau hypersurfaces, \rm
Communications in Number Theory and Physics, Vol. 1, No.3 (2007),
539-578.




\bibitem{WAN93}
D. Wan, \it Newton polygons of zeta functions and L functions, \rm
Ann. of Math, Vol. 2,  No. 2(1993),  249-293.

\bibitem{WAN99} D. Wan,
\it Dwork's conjecture on unit root zeta functions, \rm Ann.
Math., 150(1999), 867-927.

\bibitem{WAN1}D. Wan,
\it Higher rank case of Dwork's conjecture,
\rm J. Amer. Math. Soc., 13(2000), 807-852.

\bibitem{WAN2}D. Wan,
\it Rank one case of Dwork's conjecture,
\rm J. Amer. Math. Soc., 13(2000), 853-908.

\bibitem{WAN03}D. Wan,
\it Rationality of partial zeta functions,
\rm Indagationes Math., New Ser., 14(2003), 285-292.

\bibitem{WAN04} D. Wan,
\it Variations of $p$-adic Newton polygons for L-functions of
exponential sums, \rm Asian J. Math., Vol.8, 3(2004), 427-474.

\bibitem{WAN05} D. Wan,
\it Geometric moment zeta functions, \rm Geometric Aspects of
Dwork Theory, Walter de Gruyter, 2004 , Vol II, 1113-1129.

\bibitem{WAN06}
D. Wan, \it Arithmetic mirror symmetry, \rm Pure Appl. Math. Q.,
1(2005), 369-378.

\bibitem{WAN07}
D. Wan, \it Mirror symmetry for zeta functions, \rm Mirror
Symmetry V, AMS/IP Studies in Advanced Mathematics, Vol.38,
(2007), 159-184.

\bibitem{Yu}
J-D, Yu, Variation of the unit roots along the Dwork family of
Calabi-Yau varieties, preprint, 2007.

\end{thebibliography}
\end{document}